\documentclass[12pt]{article}

\usepackage{latexsym}
\usepackage{graphicx}
\usepackage{epsfig}
\usepackage{slovak}
\usepackage[cp1250]{inputenc}

\begin{document}

\setcounter{page}{62}

\begin{center}

     {\large\bf \uppercase{Identification of Fractional-Order
                           dynamical systems}}

        \vskip 1cm

      { \v{L}ubom\'{\i}r Dor\v{c}\'ak, Vladim\'{\i}r Le\v{s}ko,
        Imrich Ko\v{s}tial }

\vspace*{0.5cm}
\hspace*{0.0cm} Department of Management and Control Engineering\\
\hspace*{0.0cm}   BERG Faculty, Technical University Ko\v{s}ice\\
\hspace*{0.0cm}   Bo\v{z}eny N\v{e}mcovej 3, 042 00 Ko\v{s}ice, Slovakia\\
\hspace*{0.0cm} e-mail: Lubomir.Dorcak@tuke.sk,
                phone: (+42155) 6025172 \\

        \vskip 1cm

\end{center}

\noindent
{\bf \uppercase{ABSTRACT:}}
This contribution deals with identification of fractional-order
dynamical systems. We consider systems whose mathematical
description is a three-member differential equation in which
the orders of derivatives can be real numbers. We give a
discretization method and a numerical solution of differential
equations of this type. An experimental method of identification
is given which is based on evaluation of transfer characteristics.
This is a combination of the method of derivatives of transfer
characteristics and of the method of passive search. The
verification was performed on systems with known parameters and
also on a laboratory object.

\section*{\normalsize\bf \uppercase{1. Introduction}}

Real objects in general are fractional-order systems, although
in some types of systems the order is very close to an integer
order. Such systems are mainly electronic systems composed of
quality electronic elements. So far, however, the systems were
described as integer-order systems, regardless of the negative
consequences caused by neglecting the real order of the system [4,5].
Disregarding the fractional order of the system was caused mainly
by the nonexistence of simple mathematical tools for the description
of such systems. Since major advances have been made in this area
in the last few years [1,2,3,4,5,6] it is possible to consider also
the real order of the dynamical systems.
Such models are more adequate for the description of dynamical
systems with distributed parameters than integer-order models
with concentrated parameters. With regard to this,
in the task of identification, it is necessary to consider also
the fractional order of the dynamical system. In this contribution
we will concentrate mainly on the identification of parameters
(including the order of derivatives) for a chosen structure of the
model with emphasis on the methods of evaluation of transfer
characteristics. The verification of the correctness of parameter
identification will be done by using it in systems with known
parameters.

\section*{\normalsize\bf \uppercase{2. Definition of the
                            fractional-order system model}}

For the purposes of later development let us consider as the model
of a dynamical system the following
fractional-order three-member differential equation
\begin{equation} \label{r1}
      {a_2y^{(\alpha)}(t) + a_1y^{(\beta)} (t) + a_0 y(t)} = u(t)
\end{equation}
where
 $\alpha$  and $\beta$
are in general real, with initial conditions
 $y^{(\beta)}(0)=0$
 and
 $y(0)=0$.

\section*{\normalsize\bf \uppercase{3. Numerical solution
            of the fractional-order
\newline
\hspace*{0.4cm} differential equation}}

For the numerical calculation of the fractional-order diffrential
equation  (\ref{r1}) we use for an approximation of fractional
derivatives in the equation (\ref{r1}) the following formula taken
from [6]
\begin{equation} \label{r2}
   y^{(\alpha)}(t) \approx
   \; _{(t-L)}D^{\alpha}_{t}y(t)\ = \;
   h^{-\alpha} \sum_{j=0}^{N(t)}b_{j}y(t-j h ) \; ,
\end{equation}
where $L$ is the "memory length",  $h$ is the time step,
\begin{displaymath}
N(t) = \min \left\{
           \left[
              \frac{\mbox{$t$}}{\mbox{$h$}}
           \right] , \;
           \left[
              \frac{\mbox{$L$}}{\mbox{$h$}}
           \right]
         \right\} \; ,
\hspace*{1em}
\end{displaymath}
$[z]$ is the integer part of $z$,
\begin{equation} \label{r3}
       b_j = (-1)^j {{\alpha} \choose j}
\end{equation}
where
$
        {{\alpha} \choose j}
$
 is a binomial coefficient.

The approximation of the equation (\ref{r1}) in discrete time steps
 $t_m\  (m=2,3,...)$
has the following form [4]
\begin{equation} \label{r4}
       a_2h^{-\alpha}\sum_{j=0}^{m} {b_j y_{m-j}} +
       a_1h^{-\beta} \sum_{j=0}^{m} {c_j y_{m-j}} +
       a_0y_{m} = u_m
\end{equation}

From the approximation (\ref{r4}) we can derive [4] the following explicit
recursive formula for the calculation of the values
 $y_m\ (m=2,3,...)$
\begin{equation} \label{r5}
       y_m =
       \frac{u_m -
       a_2h^{-\alpha} \sum_{j=1}^{m} {b_j y_{m-j}} -
       a_1h^{-\beta}  \sum_{j=1}^{m} {c_j y_{m-j}}
            }
            {
       a_2h^{-\alpha} b_0 +
       a_1h^{-\beta}  c_0 + a_0
            }
\end{equation}
with~$y_0 = 0,\ y_1 = 0$.
An analytical solution of fractional-order differential equations
is described in detail in [2], and for closed regulation circuits
also in [4,5].

\section*{\normalsize\bf \uppercase{4. Method of differentiation
                                     of transfer
\newline
\hspace*{0.4cm} characteristics}}

Assume the input of the system
 (\ref{r1})
with known values $\alpha$  and  $\beta$
is acted upon by arbitrary function $u(t)$, and $y(t)$ is the output
function of the system, while the initial conditions need not be
zero in this method.
If the input and output functions are known functions of time
obtained, for example, by measurement, then for given discrete
times
 $t_m\  (m=2,3,...)$
 we can determine the values of the derivatives of the
output quantity according to the formula  (\ref{r2}).
 For each time instance $t_m$ the equation
 (\ref{r4})
must be satisfied. Therefore after substituting
the corresponding values into the equation
 (\ref{r4})
we obtain a system of three linear equations for unknown
coefficients $a_0, a_1, a_2$.
 By solving this system of equations and obtaining the
values of coefficients $a_0, a_1, a_2$ the system, from which we
used experimental data in the calculations, is identified.
In order to be able to use more measured values $y_i, u_i$, and to
make the identification more accurate, we consider for the estimate
of the vector of unknown parameters $\bar{a}$ as an error criterion
the following functional
\begin{equation} \label{r6}
      E(\bar{a}) =
                  \int \limits_{0}^T
      [{a_2y^{(\alpha)}(t) + a_1y^{(\beta)} (t) + a_0 y(t)} -
        u(t)]^2 dt
      \approx   min
\end{equation}
A necessary condition for the minimum to be achieved is
\begin{equation} \label{r7}
   \frac {\partial E(\bar{a})}{\partial \bar{a} } = 0
\end{equation}
Hence, after manipulations we obtain a system of three linear
equations for the vector of unknown parameters $\bar{a}$
\begin{eqnarray*}
\; {a_2 \int \limits_{0}^T (y^{(\alpha)}(t))^2 dt +
    a_1 \int \limits_{0}^T  y^{(\alpha)}(t) y^{(\beta)}(t) dt +
    a_0 \int \limits_{0}^T  y^{(\alpha)}(t) y(t) dt} =
             \int \limits_{0}^T y^{(\alpha)}(t) u(t) dt
\end{eqnarray*}
\begin{equation} \label{r8}
   {a_2 \int \limits_{0}^T  y^{(\alpha)}(t) y^{(\beta)}(t) dt +
    a_1 \int \limits_{0}^T (y^{(\beta)}(t))^2 dt +
    a_0 \int \limits_{0}^T  y^{(\beta)}(t) y(t) dt} =
             \int \limits_{0}^T y^{(\beta)}(t) u(t) dt
\end{equation}
\begin{eqnarray*}
   {a_2 \int \limits_{0}^T  y^{(\alpha)}(t) y(t) dt \;\; +\;\;
    a_1 \int \limits_{0}^T  y^{(\beta)}(t)  y(t) dt \;\; +\;\;
    a_0 \int \limits_{0}^T (y(t))^2 dt} \; = \;
             \int \limits_{0}^T y(t) u(t) dt
\end{eqnarray*}
where $y(t)$ and $u(t)$ are the measured output and input functions of
the system, and the derivatives $y^{(\alpha)}(t)$ and $y^{(\beta)}(t)$
are the derivatives of the output function of the system, obtained by
calculation from the formula (\ref{r2}).
In the discretization of the system (\ref{r8})  we substitute the
continuous functions
$y(t), y^{(\alpha)}(t), y^{(\beta)}(t),$ $u(t)$
in discrete time instants $t_m (m=0,1,2,3, ... ,M)$ with their
discrete values $y_m, y^{(\alpha)}_m, y^{(\beta)}_m, u_m$ and the
integration we substitute with the summation for the time interval
 $T$ for which we have $M+1$ measured values from
the system. Then we get
\begin{eqnarray*}
\; {a_2 \sum \limits_{m=0}^M (y^{(\alpha)}_m)^2  +
    a_1 \sum \limits_{m=0}^M  y^{(\alpha)}_m y^{(\beta)}_m  +
    a_0 \sum \limits_{m=0}^M  y^{(\alpha)}_m y_m } \; = \;
             \sum \limits_{m=0}^M y^{(\alpha)}_m u_m
\end{eqnarray*}
\begin{equation} \label{r9}
   {a_2 \sum \limits_{m=0}^M  y^{(\alpha)}_m y^{(\beta)}_m  \; + \;
    a_1 \sum \limits_{m=0}^M (y^{(\beta)}_m)^2  \; + \;
    a_0 \sum \limits_{m=0}^M  y^{(\beta)}_m y_m } \; = \;
             \sum \limits_{m=0}^M y^{(\beta)}_m u_m         \;\;\;\;
\end{equation}
\begin{eqnarray*}
   {a_2 \sum \limits_{m=0}^M  y^{(\alpha)}_m y_m  \;\; +\;\;
    a_1 \sum \limits_{m=0}^M  y^{(\beta)}_m  y_m  \;\; +\;\;
    a_0 \sum \limits_{m=0}^M (y_m)^2} \; = \;
             \sum \limits_{m=0}^M y_m u_m
\end{eqnarray*}

From the system of equations (\ref{r9}) we can calculate the unknown
parameters  $\bar{a}$, and the task of parameter identification of the
system is thus solved provided we know the values
 $\alpha$  and $\beta$.

\section*{\normalsize\bf \uppercase{5. Method of identification
                                       of fractional-order
\newline
\hspace*{0.4cm} dynamical systems}}

 In general
in the sytems of the type
  (\ref{r1}) it is necessary to identify the orders of derivatives
 $\alpha$ and $\beta$ which are non-integer and real. If, under this
assumption, we used the given method of differentiation of transfer
characteristics for the system
 (\ref{r1}), we would obtain a system of five nonlinear equations.
It is difficult to solve such systems of equations, however, therefore
we chose a combination of the method of differentiation of transfer
characteristics for the calculation of the coefficients
 $a_0, a_1, a_2$  and the method of simple search
for the determination of the coefficients $\alpha$ and $\beta$. While
the first method is analytical, in the second the coefficients are
calculated iteratively.

In the identification of the parameters of the model we proceed as
follows:
\begin{description}
\item[ (a)]
      For the parameters $\alpha$ and $\beta$ we determine the intervals
      $<\alpha_{min},\alpha_{max}>$ and $<\beta_{min},\beta_{max}>$,
      in which we look for the parameters. We divide each of these
      intervals into $n_\alpha$ and $n_\beta$ subintervals
   \begin{equation} \label{r10}
      n_\alpha = \frac {2 (\alpha_{max} - \alpha_{min})}
                       {\varepsilon}
   \end{equation}
   \begin{equation} \label{r11}
      n_\beta = \frac {2 (\beta_{max} - \beta_{min})}
                       {\varepsilon}
   \end{equation}
      where $\varepsilon$ denotes the fineness of division. We choose
      the starting values $\alpha$ and $\beta$ from the chosen intervals.
\item[   (b)]
      We use the method of differentiation of the transfer
      characteristic to calculate for given $\alpha$ and $\beta$
      the corresponding values of the parameters $a_0, a_1, a_2$.
\item[   (c)]
      We calculate the criterion of the approximation of experimental
      output values of the real system $ye_i$ and the outputs of the model
      $y_i$
   \begin{equation} \label{r12}
      Q = \frac {1} {M+1}
          \sum \limits_{i=0}^M [ ye_i - y_i ]^2
   \end{equation}
      which serves for the choice of optimum values of the coefficients
      $\alpha$ and $\beta$.
\item[   (d)]
      If the value of the actually computed approximation criterion
      is less than the value of the criterion from the preceding iteration
      step, we choose for the optimum subintervals those containing last
      values of the coefficients $\alpha$ and $\beta$.
\item[   (e)]
      After the calculation of the approximation criterion for  all
      subintervals, we determine from the  optimum subintervals the
      intervals
      $<\alpha_{min},\alpha_{max}>$ and $<\beta_{min},\beta_{max}>$
      and for these subintervals the new subintervals.
      We repeat this procedure until the required accuracy of the
      calculation of parameters $\alpha$ and $\beta$ is achieved, where
      we have always calculated the values of the parameters
       $a_0, a_1, a_2$ by using the method of differentiation of
      transfer characteristics.
\end{description}

\section*{\normalsize\bf \uppercase{6. Verification of the method}}

The verification of the correctness of the described method
of parameter identification will be done by using it in systems
with known parameters.

Consider \ an \ integer-order \ system [4] with \ coefficients
 \ \ $a_2=1, \  a_1=3, \ a_0=2,
\newline
 \alpha=2, \ \beta=1 .$
The transfer characteristic of this system has aperiodic behavior
and is shown in Fig.1. For the chosen intervals
$\alpha \in <2; 2>$ and $\beta  \in <1; 1>$ the precise
values of the coefficients $a_2=1,\ a_1=3,\ a_0=2$ were calculated
with the method of differentiation of transfer characteristics.
For the intervals
$\alpha \in <1,5; 2,55>$ and $\beta  \in <0,7; 1,33>$
the following values of the coefficients were calculated with both
methods $a_2=1,0001 \ ,\ $
$a_1=2,99987 ,\ a_0=1,99998 , \ \alpha=1,99993 ,
 \ \beta=0,99996$ ,
which values are very close to the parameters of the identified
system.

\begin{figure}[ht] \label{obr1_2}
    \vskip 2 mm
     \centerline{\psfig{file=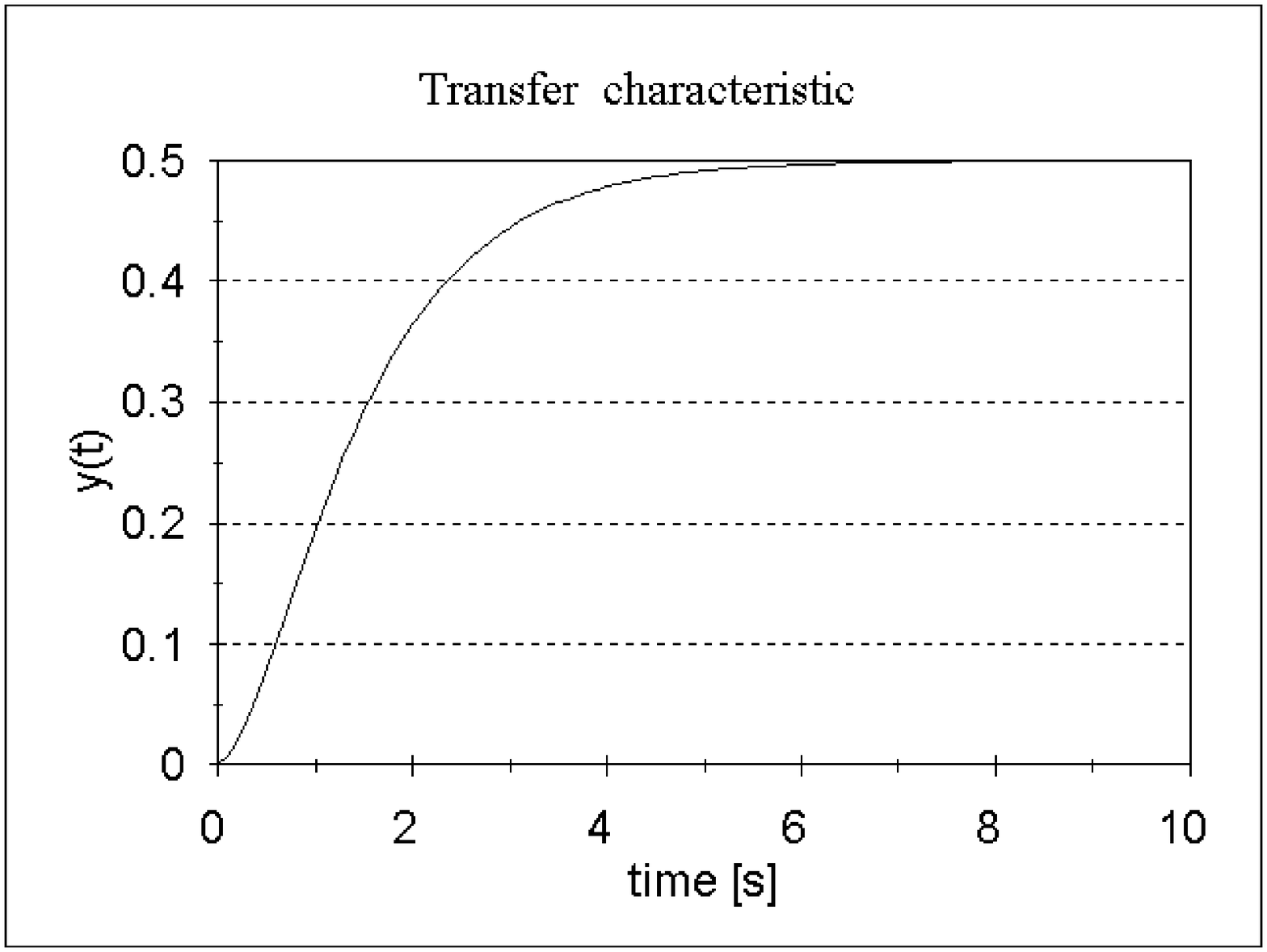,width=5.4cm} \hskip 2.2cm
		 \psfig{file=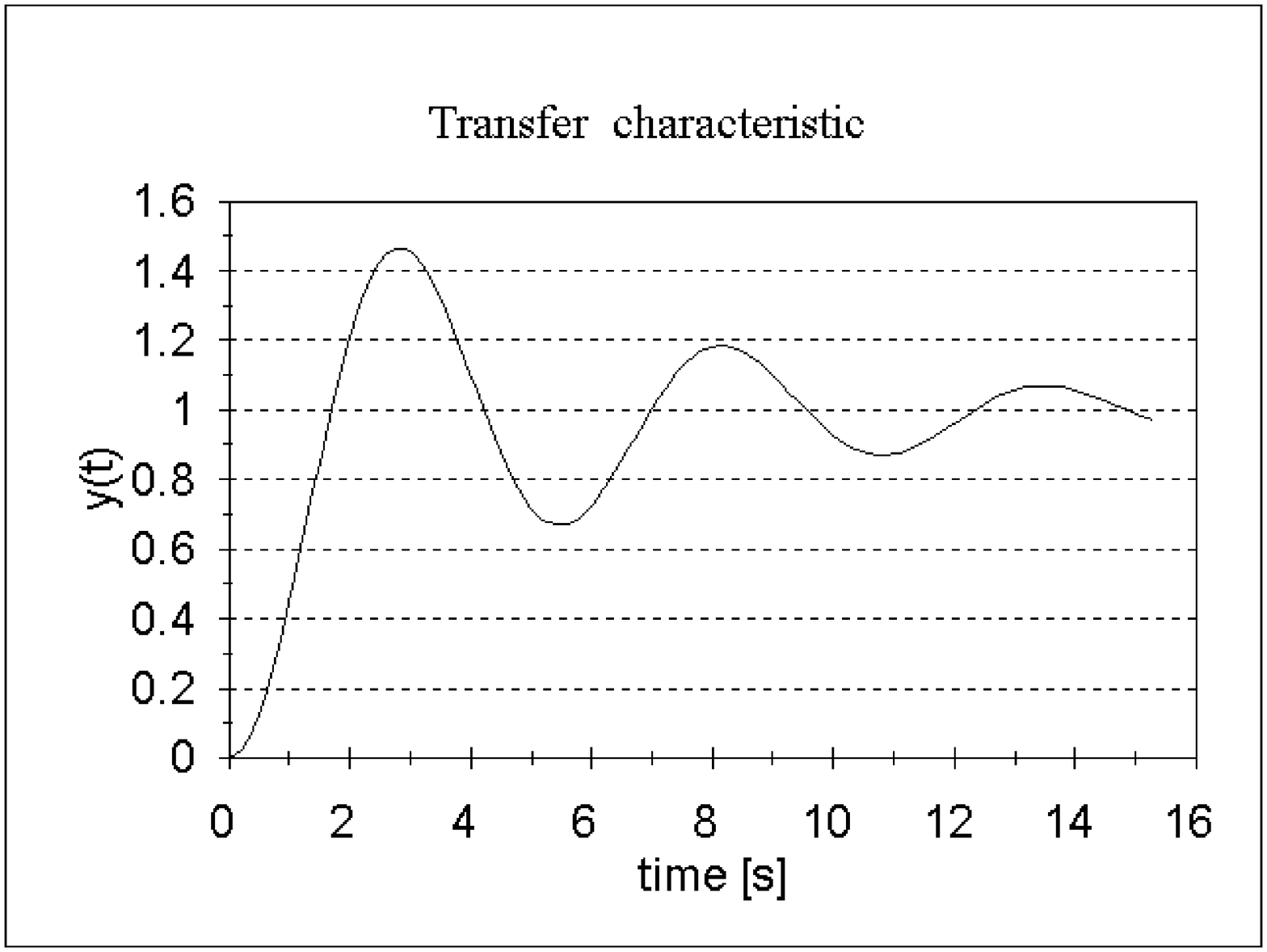,width=5.4cm}} 
    \vskip 3 mm
    \hskip 0.6cm  \rm Fig.1 Integer-order system
    \hskip 3.0cm  \rm Fig.2 Fractional-order system \rm
    \vskip 0 mm
\end{figure}

Consider \ a \ fractional-order \ system \ with \ coefficients [4] \
$a_2=0,8 \ , \ \  a_1=0,5 \ ,
\newline
 a_0=1, \ \alpha=2,2 \ , \ \beta=0,9$.
The transfer characteristic of this system computed according to
the relation (\ref{r5}),
has periodic behavior and is shown in Fig.2. For the chosen intervals
$\alpha \in <2,2; 2,2>$ and $\beta  \in <0,9; 0,9>$
the precise values of the coefficients
$a_2=0,8 \ ,\ a_1=0,5 \ ,\ a_0=1$
were computed with the method of differentiation of transfer
characteristics. For the intervals
$\alpha \in <1,5; 2,55>$ and $\beta  \in <0,7; 1,33>$
the \ following \ values \ were \ computed \ by \ using \
both \ methods
 \ $a_2=0,80005 \ ,\ $
$a_1=0,49996 \ ,\ a_0=0,99998 \ , \ \alpha=2,19996 \ ,
 \ \beta=0,89989$ \ , \
which values are very close to the parameters of the identified system.

\begin{figure}[ht] \label{obr3_4}
    \vskip 2 mm
     \centerline{\psfig{file=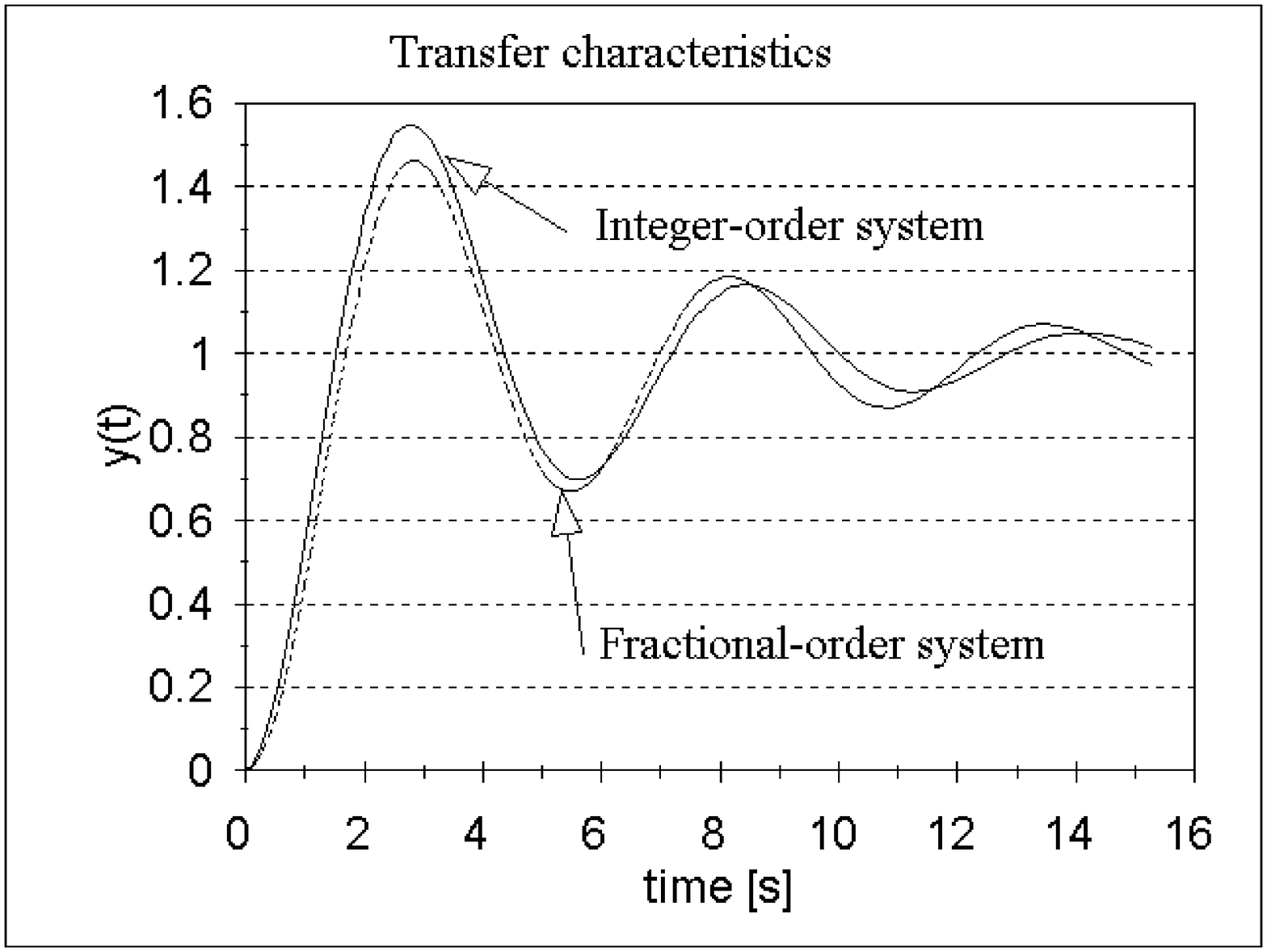,width=5.4cm} \hskip 2.2cm
		 \psfig{file=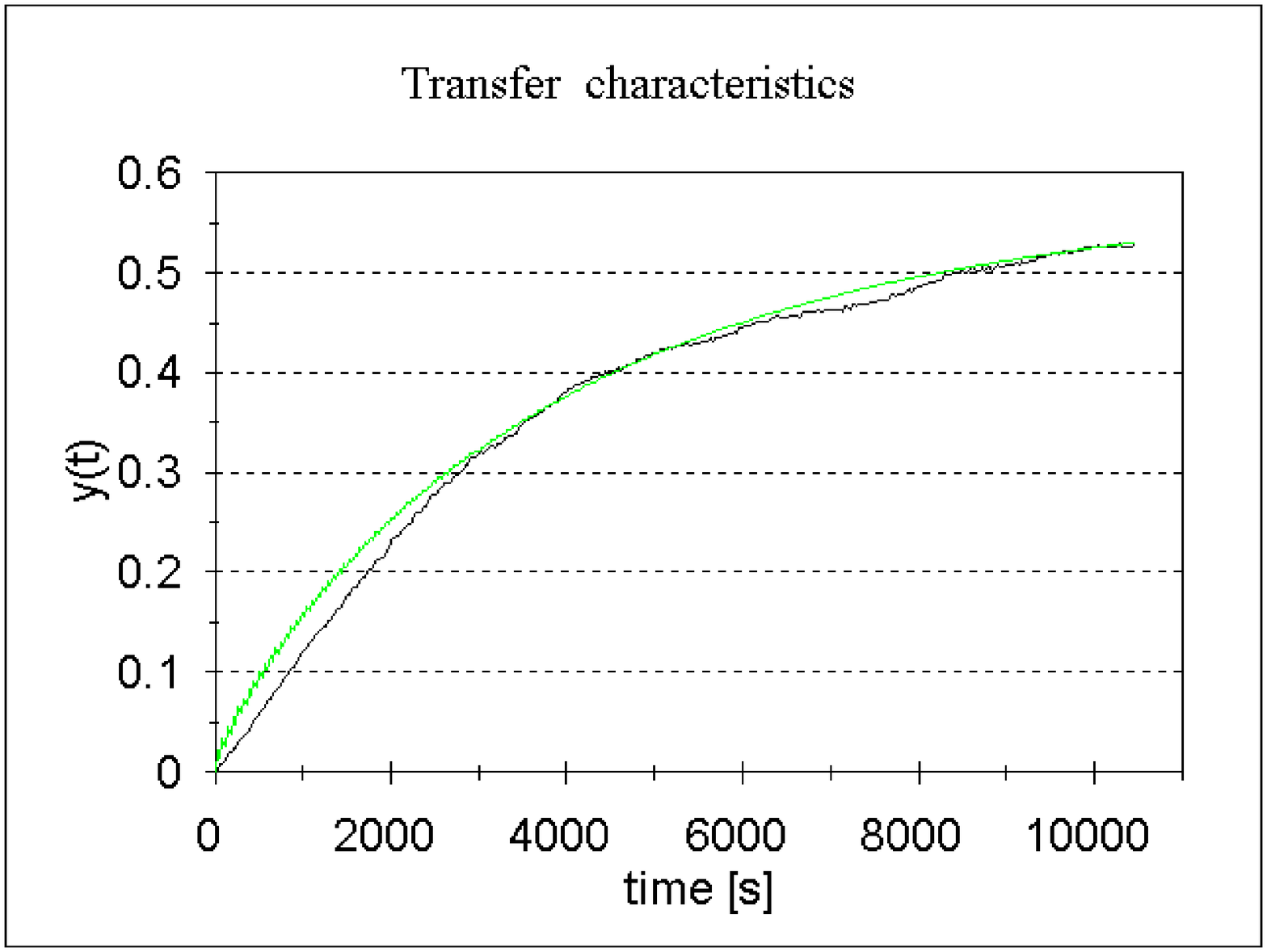,width=5.4cm}} 
    \vskip 3 mm
    \hskip 0.6cm  \rm Fig.3 Approximation of system
    \hskip 3.0cm  \rm Fig.4 Real object \rm
    \vskip 0 mm
\end{figure}

Let us consider again the previous fractional-order system with
coefficients $a_2=0,8 \ , $
$ a_1=0,5 \ , \ a_0=1, \ \alpha=2,2 \ , \ \beta=0,9 \ ,$
 whose behavior is given in Fig.2. Let us identify
it as a second-order system, which was the classical procedure
up until now. If we  want to identify the system as a second-order
system, we have to specify the same upper and lower bounds
$\alpha \in <2; 2>$ and $\beta  \in <1; 1> .$
 With the method of
differentiation of transfer characteristics the following coefficient
values were calculated :
$a_2=0,76639 \ , \ a_1=0,23184 \ , \ a_0=1 .$
The system with these parameters is actually a classical integer-order
model of a fractional-order system; their graphical comparison is
shown in Fig.3. Although optically both characteristics are rather
close to each other, the use of such model in the synthesis of the
parameters of the regulator is associated with unfavorable
consequences for the dynamics of the regulation circuit as shown
e.g. in [4,5].

Let us now identify the parameters of an experimental heating
furnace whose transfer characteristic is shown in Fig. 4.
Assuming it is a second order system with $\alpha=2, \
 \beta=1, $ we obtain parameters $ \ a_0=1,928 \ , \
a_1=4892,733 \ , \ a_2=-73043,36$ and the value $1,02.10^{-3}$ of the
criterion (\ref{r12}).
For the chosen intervals
 $\alpha \in <1,1; 2,55>$ and $\beta  \in <0,33; 1,3>$ ,
the parameters
 $a_2=-14994,3 \ ,\ $
 $a_1=6009,52 \ ,
 \ a_0=1,69 \ , \ \alpha=1,31 \ , \
  \beta=0,97$
were identified for this system, and the criterion  (\ref{r12})
had the value $2,7.10^{-4}$  (see Fig. 4).
 As we can see ($a_2$) the previous two approximations are unsuitable.
Assuming a two-member differential equation, the parameters take
on the values
 $a_1=788,35 \ ,
 \ a_0=1,39 \ , \
  \beta=0,73$
and the value of the criterion (\ref{r12}) is $6,3.10^{-4}$.

\section*{\normalsize\bf \uppercase{7. Conclusion}}

The above results of parameter identification for fractional- or
integer-order dynamical systems confirm that in known systems the
agreement of calculated and actual parameters was very good.
In the case of real objects identification it is necessary to pay
attention to the accuracy of the output quantity measurement
in the identified system, and to the time step of the measurement
from the point of view of algorithms.

When identifying the parameters of fractionarl-oder dynamical systems
in which real orders of derivatives $\alpha$ and $\beta$ are identified,
it is suitable, considering the existence of local minima of the
functional (\ref{r6}) ,
to repeat the calculations for different intervals
 $<\alpha_{min},\alpha_{max}>$ and $<\beta_{min},\beta_{max}>$,
in which we look for the parameters, and for different fineness of
division $\varepsilon$ in the formulas
 (\ref{r10}) and (\ref{r11}).

\section*{\normalsize\bf \uppercase{8. References}}

\noindent
[1] K.~B.~Oldham and J.~Spanier: {\em The Fractional Calculus}.
    Academic Press, New York, 1974.

\noindent
[2] I.~Podlubny : {\em The Laplace Transform Method for Linear
    Differential Equations of the Fractional Order},
    UEF - 02 - 94, The
    Academy of Sciences Institute  of Experimental Physics, 1994,
    Kosice.

\noindent
[3] M.~Axtell and M.~E.~Bise : {\em Fractional Calculus Applications
    in Control Systems}, in: Proc.  of the IEEE 1990 Nat. Aerospace
    and Electronics Conf., New York, 1990, pp. 563 - 566.

\noindent
[4] \v{L}. Dor\v{c}\'ak :
    {\em Numerical Models for Simulation of the Fractional-Order
         Control Systems}, in:
    UEF - 04 - 94, The Academy of Sciences Institute of Experimental
    Physics, 1994, Kosice, p. 12.

\noindent
[5] \v{L}. S\'ykorov\'a, I. Ko\v{s}tial, \v{L}. Dor\v{c}\'ak :
    {\em A comparison of Integer-
    and Fractional- Order PID Regulators}, in: Proceeding of the 2nd
    scientific-technical conference with international participation
    PROCESS CONTROL, Horn\'a Be\v{c}va, Czech Republic, June 3-6, 1996,
    pp. 347-351.

\noindent
[6] I.~Podlubny.
    {\em  Numerical  methods  of  the  fractional calculus},
    part I, II. Trans. of the TU Kosice, vol.~4, no.~3-4, pp. 200-208.

\end{document}